\documentclass[a4paper]{paper_cw}

\usepackage{helvet}
\usepackage{amsmath}
\usepackage{amssymb}
\usepackage{mathrsfs}
\usepackage{graphics}
\usepackage{natbib}
\usepackage{color}
\usepackage{fancyhdr} % must be loaded for header/footer
\usepackage{colortbl} % must be loaded for \columncolor
\usepackage{subfigure}
\usepackage[dvips]{graphicx} % must be loaded to scale graphics
\usepackage{marvosym}
\usepackage{natbib}

\usepackage{array} 

\usepackage{pgfplots}
\usepackage{framed}

\usepackage{tikz}
\usetikzlibrary{matrix,calc}

\long\def\symbolfootnote[#1]#2{\begingroup \def\thefootnote{\fnsymbol{footnote}}\footnote[#1]{#2} \endgroup} 

\renewcommand{\vec}[1]{ \ensuremath{ \mathbf{ #1 } } }
\newcommand{\ten}[1]{ \ensuremath {\mathbf{#1} } }

\newcommand{\tenD}{\ensuremath{ \ten{D} }}

\newcommand{\tenF}{\ensuremath{ \ten{F} }}

\newcommand{\tena}{\ensuremath{ \ten{a} }}
\newcommand{\tenb}{\ensuremath{ \ten{b} }}

\newcommand{\tenf}{\ensuremath{ \ten{f} }}

\newcommand{\tenn}{\ensuremath{ \ten{n} }}

\newcommand{\tenu}{\ensuremath{ \ten{u} }}
\newcommand{\tenv}{\ensuremath{ \ten{v} }}

\newcommand{\tenx}{\ensuremath{ \ten{x} }}

\newcommand{\vecn}{\ensuremath{ \vec{n} }}

\newcommand{\vecu}{\ensuremath{ \vec{u} }}
\newcommand{\vecv}{\ensuremath{ \vec{v} }}

\newcommand{\vecx}{\ensuremath{ \vec{x} }}

\newcommand{\onetwo}{\frac{1}{2}}

\renewcommand{\cos}[1]{ \text{cos}\hspace{0.0cm}\left( {#1} \right) }
\renewcommand{\sin}[1]{ \text{sin}\hspace{0.0cm}\left( {#1} \right) }

\newcommand{\tenone}{\ensuremath{ \ten{1} }}

\newcommand{\na}{\ensuremath{ _{n+1} }}

\begin{document}

\title{The Neural Particle Method - An Updated Lagrangian Physics Informed Neural Network for Computational Fluid Dynamics}

\author{Henning Wessels \and Christian Wei{\ss}enfels \and Peter Wriggers}

\institute{H. Wessels (\Letter) \at Institute of Continuum Mechanics, Leibniz University of Hannover, Appelstr. 11, 30167 Hannover, Germany
	\\\email { wessels@ikm.uni-hannover.de} 
}

\maketitle

\noindent This paper has been published in Computer Methods in Applied Mechanics and Engineering:\newline https://doi.org/10.1016/j.cma.2020.113127

\section*{Highlights}

\begin{itemize}
		\item A feed-forward neural network is used to construct a global geometric ansatz function.
		\item No special treatment of the incompressibility constraint is necessary.
		\item High order implicit Runge Kutta time integration is employed.
		\item Excellent conservation properties are demonstrated in numerical examples.
		\item {The computations remain stable even for irregularly distributed discretization points.}
\end{itemize}

\abstract{Today numerical simulation is indispensable in industrial design processes. It can replace cost and time intensive experiments and even reduce the need for prototypes. While products designed with the aid of numerical simulation undergo continuous improvement, this must also be true for numerical simulation techniques themselves. {Up to date, no general purpose numerical method is available which can accurately resolve a variety of physics ranging from fluid to solid mechanics including large deformations and free surface flow phenomena. These complex multi-physics problems occur for example in Additive Manufacturing processes.} In this sense, the recent developments in Machine Learning display promise for numerical simulation. {It has recently been shown  that instead of solving a system of equations as in standard numerical methods, a neural network can be trained solely based on initial and boundary conditions.} Neural networks are smooth, differentiable functions that can be used as a global ansatz for Partial Differential Equations (PDEs). While this idea dates back to more than 20 years ago \cite[]{Lagaris1998}, it is only recently that an approach for the solution of time dependent problems has been developed \cite[]{Raissi2019}. With the latter, implicit Runge Kutta schemes with unprecedented high order have been constructed to solve scalar-valued PDEs. We build on the aforementioned work in order to develop an Updated Lagrangian method for the solution of incompressible free surface flow subject to the inviscid Euler equations. The method is straightforward to implement and does not require any specific algorithmic treatment which is usually necessary to accurately resolve the incompressibility constraint. {Due to its meshfree character, we will name it the Neural Particle Method (NPM). It will be demonstrated that the NPM remains stable and accurate even if the location of discretization points is highly irregular.} }


\section{Introduction}\label{sec:introduction}

Today, numerical methods are well established and widely used in research and development. Many different methods have emerged to tackle a variety of problems. In Computational Fluid Dynamics (CFD), the study of flows in Eulerian formulation within enclosed domains, e.g. gas turbines, is typically performed with the Finite Volume (FVM) or the Finite Difference Method (FDM). In these methods, the fluid domain is discretized by a mesh on which the governing equations are solved. Standard mesh-based methods suffer from three shortcomings when applied to incompressible fluid flow:
\begin{enumerate}
	\item {In an Eulerian reference frame, convective terms occur and require  stabilization \cite[]{Patankar1980}}.
	\item To represent continuous free surface flows, mesh-based methods require special techniques such as the Volume-of-Fluid (VoF) approach \cite[]{Hirt1981}.
	\item The incompressibility constraint causes numerical instabilities and must be stabilized {\cite[]{Brezzi1974}}.
\end{enumerate}
In order to relax the mesh-dependency and to account for free surface flows, the Arbitrary-Lagrangian-Eulerian (ALE) formulation has been introduced by \cite{Hirt.1974}, see also \cite{Tezduyar1992} and \cite{Braess2000}. In the ALE method, the mesh is advanced in time, but independently of the fluid motion to avoid critical mesh-distortion. 

With the aim to entirely remove the mesh-dependency, so-called mesh-free or particle methods have been developed. {These methods usually employ a Lagrangian formulation.} The discretization points, referred to as particles, follow the fluid motion with the effect that convective terms disappear. Using a Lagrangian description, free surface flows are naturally represented. Beside the well-know Smoothed Particle Hydrodynamics (SPH) method which was independently introduced by \cite{Lucy.1977} and \cite{Gingold.1977}, many other schemes have emerged. Mesh-free methods do not require a fixed mesh, but also in these methods a connectivity is present. It is established by a search algorithm, which must fulfill certain topological requirements \cite[]{Liu.1997}. To equilibrate unphysical configurational forces, some methods rely on $r$-adaptivity, e.g. incompressible SPH \cite[]{Lind2012} or the Optimal Transportation Meshfree Method (OTM) suggested by \cite{Li.2010}. Alternatively, the equation of motion can be stabilized, see e.g.  \cite{Weissenfels.2018} for the OTM or \cite{Ganzenmuller.2015} for SPH. {For an overview of meshfree methods, the reader is referred to the book of \cite{Li.2007}.}

The incompressible Euler and Navier-Stokes equations solved with either mesh-based or mesh-free methods have the form of a mixed-method. Mixed methods are subject to the Ladyschenskaja-Babu{\v s}ka-Brezzi (LBB) or inf-sup condition. It requires the pressure interpolants to be of lower order than those of the velocity {\cite[]{Franca1988}}. However, in practice equal-order interpolation is often preferred. In this case, the pressure degrees of freedom lay directly on the fluid boundary. This facilitates the imposition of boundary conditions and the treatment of fluid-structure interaction. Hence, a variety of stabilization techniques have been developed which circumvent the LBB condition \cite[]{Franca1988}. According to \cite{Brezzi1997}, they can be classified  according to two different strategies. The first is to modify the bilinear form in order to achieve enhanced numerical stability without compromising consistency. The standard Galerkin method is applied to the modified equations. {Commonly known methods falling into this class are the Pressure Stabilized Petrov Galerkin method (PSPG) proposed by \cite{Tezduyar1992} and the Finite Increment Calculus (FIC) formulation \cite[]{Onate2001}.} The second strategy is to enrich the standard Galerkin method with special functions. More precisely, the space of functions is enlarged such that a comparatively coarse mesh is able to deal with the effects of unresolvable scales \cite[]{Brezzi1997}. This approach is followed in the Subgrid Scale Method, in Residual-free Bubbles or the Variational Multiscale Method, see e.g. \cite{Brezzi1997} and \cite{Hughes1998}.

{In addition to the aforementioned numerical methods, data-driven approaches powered by machine learning are increasingly finding their way into CFD. One field of application is the reconstruction of flow fields from data. The reconstruction from discrete pressure or velocity measurements can be challenging if the number of available sensors is limited \cite[]{Erichson2019}. \cite{Raissi2020} have developed an algorithm which learns the velocity and pressure fields from continuous flow visualizations such as particle image velocimetry. Data-driven techniques are also applied on the simulation of fluid flow. For example, neural networks have been developed that accurately resolve the turbulent Reynolds Averaged Navier Stokes equations. A brief overview can be found in \cite{Kutz2017}. In the context of machine learning of turbulent flows, a combination of feature extraction, i.e. Model Order Reduction (MOR) with Recurrent Neural Networks (RNNs) is commonly employed.  MOR techniques to compute the lower dimensional latent space are for example Proper Order Decomposition (POD) \cite[]{Chatterjee2000} or autoencoder \cite[]{Hinton2006}.  RNNs are based on sequence learning and therefore especially suited to learn the latent space dynamics. However, RNNs are difficult to train in practice. A commonly applied and stable variant of RNNs are Long Short Term Memory networks (LSTMs) \cite[]{Hochreiter1997}. Another alternative to standard RNNs are Neural Ordinary Differential Equations (NODE) which incorporate history implicitly and can learn from temporally scattered data \cite[]{Chen2018}. The performance of NODE and LSTM in combination with POD applied to Burgers turbulence has been studied in \cite{Maulik2020}. A detailed overview of machine learning for fluid dynamics can be found in \cite{Brunton2020}.}
	
{A remedy of the aforementioned data-driven simulation concepts is that they are physics agnostic. In order to increase the fidelity of machine learning for predicitive simulations, current research focuses on the implementation of physical constraints into formerly purely data-driven approaches. For example, \cite{Mohan2020} embed the incompressibility constraint into a convolutional autoencoder using the equivalence of convolutional neural network kernels with Finite Volume stencils. While the approach steps in the direction of Physics Informed Data-Driven Simulation, it has two severe limitations. First, stencils require a fixed lattice. This makes the application of the method to scenarios like free surface flow and large deformations cumbersome. Second, only the time independent continuity equation is included as a physical constraint. The solution of the momentum equation remains a black-box.}

{Besides their application in data-driven simulation, neural networks can also be used for the numerical solution of Partial Differential Equation (PDEs) in the absence of data. An innovative approach is to use a feed-forward neural network as a global ansatz function of a PDE.} The universal approximation theorem states that a feed-forward neural network with mild assumptions on the activation function can approximate any function \cite[]{Cybenko1989, Leshno1993}. When a neural network is designed to take the solution of a PDE as output, it acts as ansatz function of the solution. The idea has first been presented by \cite{Lagaris1998} and can be regarded as a collocation method. Note that a neural network is a continuous and smooth function which can not only be differentiated with respect to the network parameters (i.e. the weights and the biases), but also with respect to its input. If the input is the spatial position, derivatives of the network output with respect to its input are then spatial derivative operators as occurring in the PDE. These derivatives can be computed with ease even for complex network architectures using Automatic Differentiation (AD), see e.g. \cite{Griewank.2008} and \cite{Korelc.2016b}.

While the first attempt of \cite{Lagaris1998} was only applicable on rectangular domains, the method has been extended to more complex geometries \cite[]{Lagaris2000}. Thanks to the computational power increase since this pioneering work and the recent developments in Machine Learning and Neural Networks, attention in the topic has increased again. A novel approach to the solution of PDEs on complex domains was presented by \cite{Berg2018}, who introduced a projection of the actual solution in order to satisfy boundary conditions exactly. For the solution of quasi-static mechanical problems, \cite{Samaniego2020} have used a neural network as ansatz for the displacement and formulated a loss function that minimizes the elastic energy. The method has been extended from small to large strain problems by \cite{Nguyen2020}. In all the aforementioned contributions only static PDEs have been considered. For time dependent problems, \cite{Raissi2019} suggested to exploit the structure of Implicit Runge Kutta (IRK) time integrators. The IRK stages and the final solution are considered as output neurons of the neural network.

Building on the aforementioned pioneering contributions, in this work the Neural Particle Method (NPM) is presented as an innovative approach to {solve} the incompressible, inviscid Euler equations. In order to account for large deformations and free surface flow, the Euler equations are formulated in an  Updated Lagrangian frame.  Velocity and pressure are approximated with a neural network as a global ansatz function. Following \cite{Berg2018}, the boundary conditions can be fulfilled exactly. The temporal discretization is realized with the high order IRK integration scheme of \cite{Raissi2019}. It will be demonstrated that the NPM fulfills the incompressibility constraint without any stabilization even on arbitrarily distributed discretization points. This is a severe benefit of the NPM compared to state of the art numerical methods.

The outline of this paper is as follows: First, the IRK method of \cite{Raissi2019} is reviewed for second order Ordinary Differential Equations (ODEs). The NPM is developed step by step in section \ref{sec:UL}. In section \ref{sec:results}, the method is validated and its performance illustrated by means of numerical examples of practical relevance, namely sloshing in a container and the classical dam break test case. The paper concludes in section \ref{sec:conclusion}.

\section{A Physics Informed Neural Network for Second Order ODEs}\label{sec:IRK_ODE}

\begin{figure}
	\centering
	\includegraphics[width=0.8\textwidth]{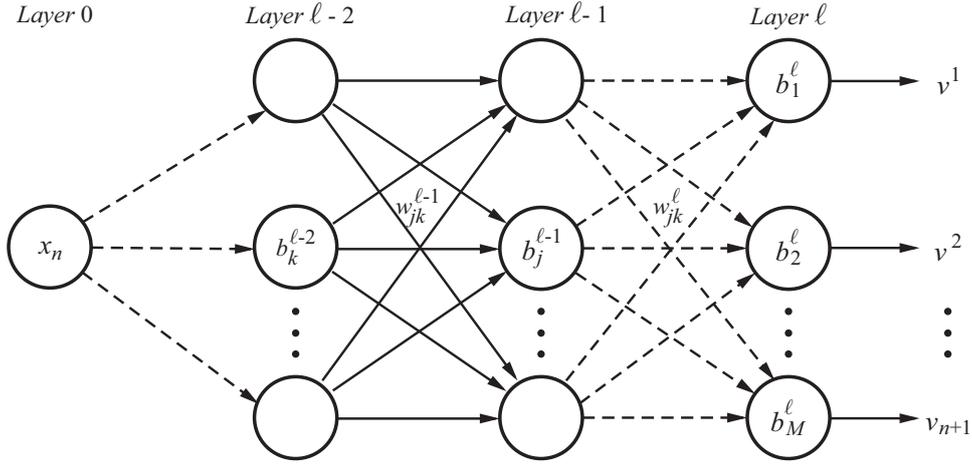}
	\caption{Schematic of a one dimensional Physics Informed Neural Network adapted from \cite{Berg2018}. The network inputs the known position $x_n$ and outputs the current velocity $v\na$ as well as the IRK velocity stages $v^j$.}
	\label{fig:FFNN}
\end{figure}

Computational mechanics is about the {differential} formulation and solution of the equation of motion, which is a second order PDE in time. In the presence of damping, the latter can be expressed in terms of the acceleration $\tena$ which is a function of time $t$, velocity $\tenv$ and  position $\vecx$:
\begin{equation}\label{eq:2ndPDE}
\begin{aligned}
\tena = \tena\left(\tenv\left(\tenx\left(t\right)\right), \tenx\left(t\right), t\right)
\end{aligned}
\end{equation}
The idea of Runge Kutta integration is to evaluate the acceleration $\tena$ at a distinct number of stages $s$ between two time steps $t_n$ and $t\na$, i.e. at time instances $t_n+c^i\Delta t$. The coefficients $c^i$ are given by the Butcher tableau along with the coefficients $a^{ji}$ and $b^j$ appearing in (\ref{eq:velocityRK}) and (\ref{eq:positionRK}). The indices $i$ and $j$ of the coefficient denote the Runge Kutta stage and range from one to the total number of stages $s$.
Implicit Runge Kutta (IRK) methods {are stable even when applied to stiff equations.} Higher orders of accuracy can be achieved by increasing the number of IRK stages $s$. Details on the theoretical background can be found in \cite{Iserles2012} or \cite{Hairer2010}. For the sake of brevity, the following abbreviations for the position and velocity stages are introduced:
\begin{equation}\label{eq:abbr_stages}
\begin{aligned}
	\tenx^i = \vecx\left(t_n + c^i \Delta t\right),
	\qquad 
	\tenv^i = \vecv\left( \vecx\left(t_n + c^i \Delta t\right), t_n + c^i \Delta t \right)
\end{aligned}
\end{equation} 
Additionally, the time dependency of a variable will be expressed by a subscript, e.g. $\vecv_n=\vecv\left(t_n\right)$. The velocity stages and the current velocity then follow the update formula:
\begin{equation}\label{eq:velocityRK}
\begin{aligned}
\tenv^j &= \vecv_n + \Delta t \sum_{i=1}^s a^{ji} \, \tena^i\left(\tenv^i, \, \tenx^i, \, t_n + c^i \Delta t\right) \\
\tenv\na &= \tenv_n + \Delta t \sum_{j=1}^s b^j \tena^j \left(\tenv^j, \tenx^j,  \, t_n + c^i \Delta t\right)
\end{aligned}
\end{equation}
{The update of the position is directly computed from the velocity stages:}
\begin{equation}\label{eq:positionRK}
\begin{aligned}
\tenx^j &= \vecx_n + \Delta t \sum_{i=1}^s a^{ji} \, \tenv^i\\
\tenx\na &= \tenx_n + \Delta t \sum_{j=1}^s b^j \tenv^j \\
\end{aligned}
\end{equation}
The update formulae (\ref{eq:velocityRK}) and (\ref{eq:positionRK}) yield a fully coupled system of equations of the size of the number of IRK stages $s \times s$. It can be solved e.g. by means of the Newton-Raphson method. For standard mesh-based and mesh-free methods this is usually too expensive. Hence, lower order time integrators such as the second order Newmark scheme are preferred. {When} a neural network is used as global spatial ansatz, \cite{Raissi2019} suggested to put a neural network prior on the Runge Kutta stages and the solution. This allows the use of very high order IRK methods.  For the IRK integration of the equation of motion (\ref{eq:2ndPDE}), the neural network takes the position of a discretization point $\vecx_n$ as input and outputs its IRK velocity stages $\tenv^i$ as well as its velocity $\tenv\na$ at the next time step. Note that the position is solely determined by the velocity through (\ref{eq:positionRK}) and does not require a separate network output. The concept is briefly sketched in figure \ref{fig:FFNN}. The parameters of the neural network are its weights $w_{jk}^l$ and its biases $b^l_j$. In a Feed-Forward Neural Network (FFNN), information is passed in one direction from the input towards the output. The feed-forward algorithm for computing the  output $v^l$ for a given input $x_n$ is defined by:
\begin{equation}
\begin{aligned}
v^l &= \sigma^l \left(z^l\right) \\
z^l &= W^l \, \sigma^{l-1} \left(z^{l-1}\right) + b^l \\
z^{l-1} &= W^{l-1} \, \sigma^{l-2} \left(z^{l-2}\right) + b^l \\
&\, \, \,\vdots\\
z^2 &= W^2 \, \sigma^1\left(z^1\right) + b^2 \\
z^1 &= W^1 \, x_n + b^1
\end{aligned}
\end{equation}
Here, $\sigma^l$ denotes the activation function which is the hyberbolic tangent throughout this work. The components of the weight matrices $W^l$ are $w^l_{jk}$. Further details can be found e.g. in \cite{Berg2018}. After initialization, the network parameters, i.e. the weights and the biases are optimized in a training loop to reduce a predefined objective loss function. In the case of Physics Informed Neural Networks, the loss function is designed to incorporate the governing physics. For this purpose, the velocity update (\ref{eq:velocityRK}) is rearranged such that it yields an expression {for the known velocity at the previous time step}: 
\begin{equation}\label{eq:approx_vn}
\begin{aligned}
\vecv_{n}^j &= \tenv^j - \Delta t \sum_{i=1}^s a^{ji} \, \tena^i \left(\tenv^i, \tenx^i, t_n + c^i \Delta t \right) \\
\tenv_n^{s+1} &= \tenv\na - \Delta t \sum_{j=1}^s b^j \tena^j \left(\tenv^j, \tenx^j, t_n + c^i \Delta t \right)
\end{aligned}
\end{equation}  
{The velocity estimates $\vecv_n^j$ and $\vecv_n^{s+1}$ are computed from the velocity IRK stages and the final solution, respectively. Therefore, the estimates are dependent on the network output. To optimize the latter, a loss function {$SSE_\vecv$} is introduced as the sum of squared errors of the velocity estimates:}
\begin{equation}\label{eq:loss}
\begin{aligned}
{SSE_\vecv = \sum_{j=1}^{s}\left|\vecv_n^j - \vecv_n \right|^2 + \left|\vecv_n^{s+1} - \vecv_n \right|^2}
\end{aligned}
\end{equation}
The equation of motion is then solved by training the neural network, i.e. by adjusting its weights and its biases  such that the loss (\ref{eq:loss}) reaches a minimum. Throughout this work, a combination of Adam Optimization \cite[]{Kingma2014} and L-BFGS-B gradient descent \cite[]{Liu1989} is employed for training. 
\begin{figure}
	\centering
	\includegraphics[width=0.8\textwidth]{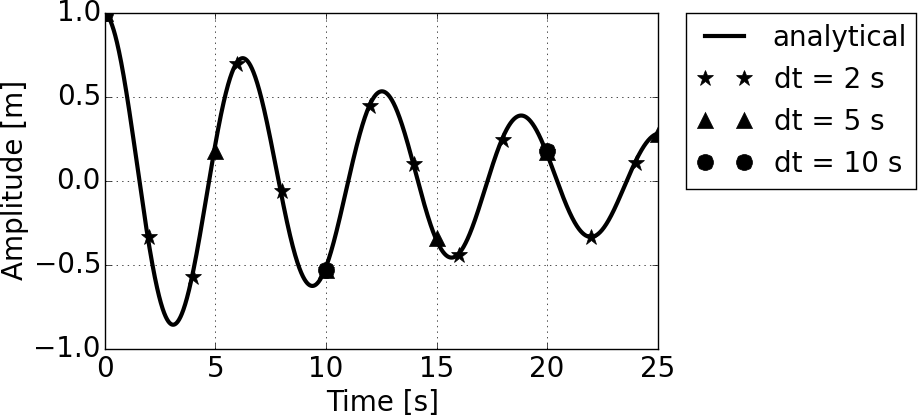}
	\caption{Comparison of Implicit Runge Kutta integration and analytical solution of the 1D equation of motion (\ref{eq:vibration_1D}) with $m=1$ kg, $k=1$ N/m and $d=0.1$ N/s. For all simulations, the network layout consists of one input neuron with 2 hidden layers of 20 neurons each and 9 output neurons, i.e $s=8$ IRK stages. With this setup, large time steps that span entire oscillations cycles can be accurately realized.}
	\label{fig:motion_1D}
\end{figure}

To illustrate the potential and the accuracy of the presented IRK integration, the initial value problem of a one dimensional mass-spring-damper system with mass $m$, stiffness $k$ and damping coefficient $d$ is considered. Written in generalized coordinates $q$, the equation of motion and its analytical solution $q\left(t\right)$ for the case of under-critical damping are given by
\begin{equation}\label{eq:vibration_1D}
\begin{aligned}
\ddot{q} + 2D\omega_0 \dot{q} + \omega_0^2 q = 0,
\qquad 
q\left(t\right) = \hat{q} \,e^{-D\omega_0 t} \, \cos{\omega_0\sqrt{1- D^2} \, t}
\end{aligned}
\end{equation}
with $\hat{q}$ the initial amplitude at time $t=0$. The eigen frequency $\omega_0$ and the attenuation factor $D$ are defined as:
\begin{equation}
\begin{aligned}
\omega_0^2 = \frac{k}{m}, \quad 
D = \frac{d}{2\sqrt{cm}}
\end{aligned}
\end{equation}
For a detailed derivation, the reader is referred to standard text books, e.g. \cite{Magnus2013}. In figure \ref{fig:motion_1D}, the IRK integrated equation of motion (\ref{eq:vibration_1D}) is plotted together with its analytical solution as a function of time. The system considered has unit mass, unit stiffness and the damping coefficient is $d=0.1$ N/s. A total number of $s=8$ IRK stages and two hidden layers with 20 neurons each have been used in all simulations. For training, 100 iterations with an Adam Optimizer (learning rate 0.001) were followed by L-BFGS-B gradient descent until convergence.  Figure \ref{fig:motion_1D} neatly illustrates the potential of high order IRK integration for computational mechanics. Time steps can be chosen so large that entire vibration cycles are spanned without losing accuracy. {This is in sharp contrast to explicit time integration schemes that are usually employed in meshfree particle methods.}

\section{Updated Lagrangian Formulation}\label{sec:UL}

So far, only the solution of an ODE has been discussed. The mass-spring-damper system from the previous section was considered as rigid. In CFD however, one is interested in the deformation of fluids which brings into play spatial derivative operators. {Inviscid incompressible fluids are described by the Euler equations.} The latter comprise the momentum equation and the balance of mass which requires the divergence of the velocity to vanish:
\begin{equation}\label{eq:EulerEquations}
\begin{aligned}
	\tena &= - \frac{1}{\rho} \operatorname{grad} p + \tenb  \\
	\operatorname{div} \vecv &= 0
\end{aligned}
\end{equation}
Here, $\rho$ denotes the {(constant)} density, $p$ the pressure and $\tenb$ the gravitational acceleration. The structure of the Euler equations imposes difficulties for its numerical solution with traditional numerical methods. To avoid instabilities arising from the incompressibility constraint, the original equations are often modified using stabilizing terms as briefly discussed in section \ref{sec:introduction}. In contrast, it will be demonstrated that the neural network architecture presented in this work is able to {compute} the pressure that fulfills the incompressibility constraint exactly without any additional algorithmic treatment. 

In order to evolve (\ref{eq:EulerEquations}) in time, a neural network is constructed that inputs analogously to section \ref{sec:IRK_ODE} the position $\vecx_n$ at time $t_n$. For each position $\vecx_n$, the network is trained to predict the IRK velocity stages $\vecv^i$, the velocity of the next time step $\vecv\na$ and the pressure stages $p^i$. Due to this network architecture, spatial derivatives can only be constructed with respect to the previous time step $t_n$. Therefore, an Updated Lagrangian formulation is employed in which the configuration at time $t_n$ is considered as the reference configuration. Spatial derivatives are computed {with respect to} the reference configuration and are then transformed into the current one via a push-forward operation. The mapping between both configurations is defined by the incremental deformation gradient $\Delta \tenF$: 
\begin{equation}
\begin{aligned}
	\frac{\partial \bullet}{\partial \vecx\na} \approx \frac{\partial \bullet}{\partial \vecx_n} \cdot \frac{\partial \vecx_n}{\partial \vecx\na} = \frac{\partial \bullet}{\partial \vecx_n} \cdot \left(\Delta \tenF\na\right)^{-1},
	\qquad
	\Delta \tenF\na = \frac{\partial \vecx\na}{\partial \vecx_n}	
\end{aligned}
\end{equation}
Note that this procedure corresponds to the chain rule of differentiation. In order to apply the IRK integration introduced in the previous section, the Euler equations must be evaluated at each IRK stage. This requires the computation of an incremental deformation gradient from each position stage $\tenx^i$. Since the neural network outputs velocity stages $\tenv^i$, first the rate of the incremental deformation gradient $\Delta \dot{\tenF}$ is computed. Making use of (\ref{eq:positionRK}), the incremental deformation gradient itself can then be obtained from its rate through IRK integration:
\begin{equation}\label{eq:defo_grad}
\begin{aligned}
&\begin{split}
\Delta \dot{\tenF}^i &= \frac{\partial \tenv^i}{\partial \vecx_n}, 
\end{split}
\qquad
&\begin{split}  
\Delta \tenF^i &= \frac{\partial \tenx^i}{\partial \vecx_n} = \tenone + \Delta t \sum_{i=1}^s a^{ji} \, \frac{\partial \tenv^i}{\partial \vecx_n} = \tenone + \Delta t \sum_{i=1}^s a^{ji} \, \Delta \dot{\tenF}^i
\end{split}\\
&\begin{split}
	\Delta \dot{\tenF}\na = \frac{\partial \tenv\na}{\partial \vecx_n}, 
\end{split}
\qquad
&\begin{split}
	\Delta \tenF\na = \frac{\partial \tenx\na}{\partial \vecx_n} = \tenone + \Delta t \sum_{j=1}^s b^{j} \, \frac{\partial \tenv^j}{\partial \vecx_n} = \tenone + \Delta t \sum_{j=1}^s b^{j} \, \Delta \dot{\tenF}^j
\end{split}
\end{aligned}
\end{equation}
Using the above defined strain measures, the velocity divergence at each IRK stage and at the next time step is computed from:
\begin{equation}\label{eq:div_v}
\begin{aligned}
	\operatorname{div} \tenv^i 
	&= \operatorname{tr} \frac{\partial \tenv^i}{\partial \tenx^i} 
	= \operatorname{tr}\left[\Delta \dot{\tenF}^i \cdot \left(\Delta \tenF^i\right)^{-1}\right] \\
	\operatorname{div} \tenv\na &= \operatorname{tr} \frac{\partial \tenv\na}{\partial \tenx\na}  = \operatorname{tr}\left[\Delta \dot{\tenF}\na \cdot \left(\Delta \tenF\na\right)^{-1}\right]
\end{aligned}
\end{equation}
The gradient of the IRK pressure stages $p^i$ is computed from:
\begin{equation}\label{eq:grad_p}
\begin{aligned}
%	\operatorname{grad} p\na &= \frac{\partial p\na}{\partial \tenx\na} = \frac{\partial p\na}{\partial \tenx_n} \frac{\partial \tenx_n}{\partial \tenx\na} =\frac{\partial p\na}{\partial \tenx_n} \, \Delta \tenF^{-1}\\
\operatorname{grad} p^i &= \frac{\partial p^i}{\partial \tenx^i} 
%= \frac{\partial p^i}{\partial \tenx_n} \frac{\partial \tenx_n}{\partial \tenx^i} 
=\frac{\partial p^i}{\partial \tenx_n} \cdot \left(\Delta \tenF\right)^{-1}
\end{aligned}
\end{equation}
An evolution equation for the pressure does not exist. While one could alternatively derive a pressure Poisson equation from the divergence of the momentum equation \cite[]{Ferziger.2002}, this was not found {to be} necessary within the NPM. The pressure is rather computed for each IRK stage such that the incompressibility constraint is fulfilled. Therefore, only the pressure stages are considered as output neurons while the current pressure is defined as the weighted average of pressure stages $p^j$:
\begin{equation}
\begin{aligned}
	p\na = \sum_{j=1}^s b^j p^j, \qquad \text{with} \quad \sum_{j=1}^s b^j = 1
\end{aligned}
\end{equation}
The current velocity $\vecv\na$ is obtained directly from the output of the neural network while the update of spatial coordinates follows (\ref{eq:positionRK}). In order to account for the mass equation, the velocity divergence is added to the loss function (\ref{eq:loss}) as the sum of squared errors:
\begin{equation}\label{eq:loss_div}
\begin{aligned}
	SSE_{\operatorname{div}\vecv} = \sum_{i=1}^{s}\left|\operatorname{div} \vecv^i\right|^2 + \left|\operatorname{div} \vecv\na\right|^2
\end{aligned}
\end{equation}
Following the original paper of \cite{Raissi2019}, Dirichlet type boundary conditions for the velocity $\tenv\left(\bar{\tenx}_\tenv\right) = \bar{\tenv}$ and pressure $p\left(\bar{\tenx}_p\right) = \bar{p}$ must also be accounted for in the loss function. Their contributions from the sum of squared errors are:
\begin{equation}\label{eq:loss_boun}
\begin{aligned}
	SSE_{\bar{\vecv}} = \sum_{i=1}^{s} \left|\vecv^i\left(\bar{\tenx}_\tenv\right) - \bar{\vecv}\right|^2 + \left|\vecv\na\left(\bar{\tenx}_\tenv\right) - \bar{\vecv}\right|^2,
	\qquad 
	SSE_{\bar{p}} = \sum_{i=1}^{s}\left|p^i\left(\bar{\tenx}_p\right) - \bar{p}\right|^2
\end{aligned}
\end{equation}
{Note that the contributions $SSE_\vecv$ (\ref{eq:loss}) and $SSE_{\operatorname{div}\vecv}$ (\ref{eq:loss_div}) involve a summation over all points and the boundary contributions (\ref{eq:loss_boun}) a summation over all boundary points. These have been omitted for brevity. The original treatment of boundary conditions suggested by \cite{Raissi2019} according to (\ref{eq:loss_boun}) is critically discussed in the remainder of this section. In addition, an approach to account for contact with rigid walls is presented.}

\subsection{Imposition of Dirichlet boundaries}\label{sec:boundary}

\begin{figure}
	\centering
	\subfigure{\includegraphics[width=0.49\textwidth]{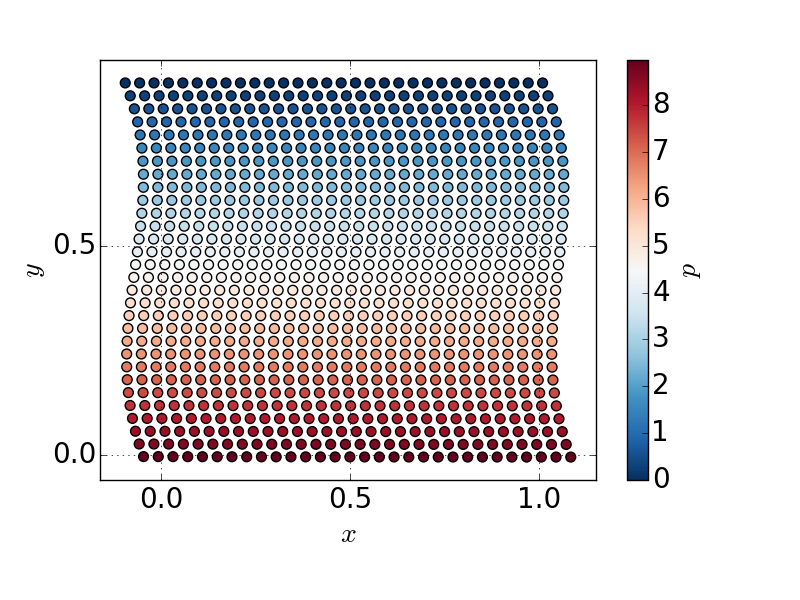}}
	\subfigure{\includegraphics[width=0.49\textwidth]{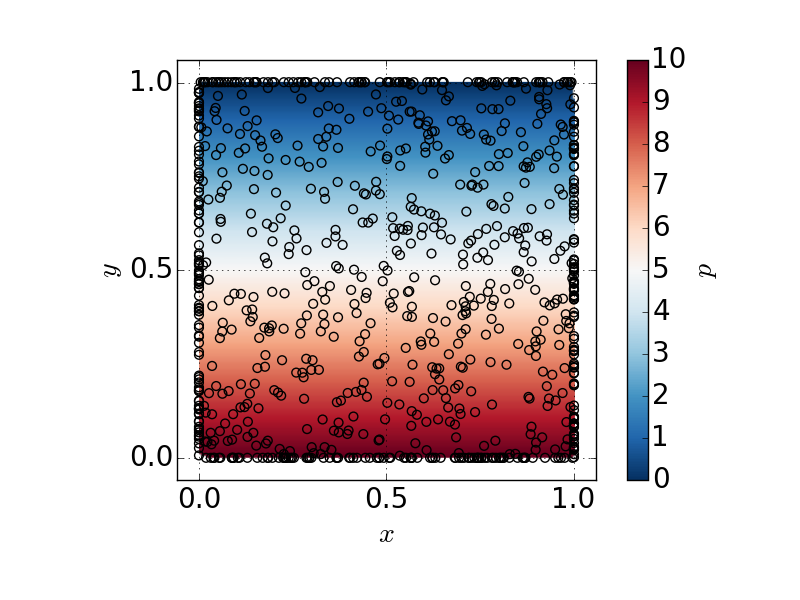}}
\caption{Static pressure of an incompressible fluid in a container after 50 time steps ($T=50$s, $\Delta t=1$s) discretized by 900 fluid particles. Left: The Dirichlet boundary conditions are not fulfilled when simultaneously trained with the incompressible Euler equations. Right: Using a velocity projection introduced by \cite{Berg2018}, the problem is solved exactly even on arbitrarily distributed discretization points (empty circles). The pressure is projected onto a background mesh using the linear interpolation method  \texttt{scipy.interpolate.griddata()} in python 3.}
\label{fig:BC_training}
\end{figure}

\cite{Raissi2019} suggest the simultaneous training of boundary data and PDE. However, when simulating multiple time steps, the accuracy of this approach is not sufficient. {This can be illustrated by a simple example. A container of unit width and height filled with an incompressible fluid of unit density under constant gravity acceleration $g=10$ m/s$^2$ is loaded only by its own weight}. As Dirichlet conditions, the velocity {normal} to the side and bottom walls must be zero as well as the pressure on the top of the domain. An accurate dynamic simulation must maintain both a linear static pressure field and the boundary conditions over time. As shown in figure \ref{fig:BC_training}, the original method of \cite{Raissi2019} fails to meet this requirement. {The fluid leaks at the boundaries.} In \cite{Berg2018}, an alternative formulation for imposing boundary conditions is suggested. The ansatz for a scalar valued {primary variable} $u\left(\tenx\right)$ can be written as the sum of a smooth extension of the boundary data $G\left(\tenx\right)$ and a smooth distance function $D\left(\tenx\right)$ multiplied with the output of a neural network $\hat{u}\left(\tenx\right)$:
\begin{equation}
\begin{aligned}
	u\left(\tenx\right) = G\left(\tenx\right) + D\left(\tenx\right) \hat{u}\left(\tenx\right)
\end{aligned}
\end{equation}
Any loss function is then formulated in terms of the ansatz $u\left(\tenx\right)$ that fulfills the boundary conditions by definition. Training of $\hat{u}\left(\tenx\right)$ is then only necessary to compute the PDE inside the domain. If for complex geometries $D\left(\tenx\right)$ and $G\left(\tenx\right)$ are difficult to define analytically, both functions can be computed using low-capacity ANNs, which ensures smoothness and differentiability. For simple rectangular geometries an analytical expression is preferred in order to save computational resources. In the static pressure test case, the normal velocity at the walls is zero and therefore the boundary extension $G\left(\tenx\right)$ vanishes. The velocity projection simplifies to
\begin{equation}\label{eq:projection}
\begin{aligned}
\begin{pmatrix}
v_x^i\left(\tenx^i\right) & v_x\left(\tenx\na\right) \\
v_y^i\left(\tenx^i\right)  & v_y\left(\tenx\na\right)
\end{pmatrix}
= \begin{pmatrix}
D_{v_x}\left(x_n \right) & D_{v_x}\left(x_n \right) \\
D_{v_y}\left(y_n \right) & D_{v_y}\left(y_n \right)
\end{pmatrix}
\circ
\begin{pmatrix}
\hat{v}_x^i\left(\tenx^i\right) & \hat{v}_x\left(\tenx\na\right)\\
\hat{v}_y^i\left(\tenx^i\right) & \hat{v}_y\left(\tenx\na\right)
\end{pmatrix}
\end{aligned}
\end{equation}
where the $\circ$ denotes the Hadamard product, i.e. element wise multiplication. Since only first order spatial derivatives appear in the Euler equations (\ref{eq:EulerEquations}), it is sufficient that the distance function is $C_1$ continuous. With $w$ the width and $h$ the height of the container, the analytical distance functions $D_{v_x}\left(x_n\right)$ and $D_{v_y}\left(y_n\right)$ introduced above for the static pressure test case can be defined as:
\begin{equation}\label{eq:velo_proj}
\begin{aligned}
D_{v_x}\left(x_n\right) = - \frac{4}{w^2}x_n^2 + \frac{4}{w} x_n,
\qquad
D_{v_y}\left(y_n\right) = \frac{1}{h}y_n
\end{aligned}
\end{equation}
Note that herein the distance functions $D_{v_x}$ and $D_{v_y}$ are computed in terms of the position $\tenx_n$ at time $t_n$, since the a-priori computation of distance functions depending on the updated position is not possible. It will be demonstrated that this approach does not alter the solution at $t\na$. As a result of the projection (\ref{eq:projection}), the rate of the incremental deformation gradient $\Delta \dot{\tenF}$ is now subject to the product rule of differentiation:
\begin{equation}\label{eq:rate_Berg}
\begin{aligned}
	\Delta \dot{\tenF}^i &= \frac{\partial \tenv^i}{\partial \vecx_n} = 
	\frac{\partial \tenD}{\partial \tenx_n}
\circ
\begin{pmatrix}
	v^i_x & v^i_x \\
	v^i_y & v^i_y
\end{pmatrix}
+ 
\begin{pmatrix}
	D_{v_x} & D_{v_x} \\
	D_{v_y} & D_{v_y}
\end{pmatrix}
\circ
\frac{\partial \hat{\vecv}^i}{\partial \vecx_n},
\qquad
\tenD = \begin{pmatrix}
D_{v_x}\left(x_n\right) \\
D_{v_y}\left(y_n\right)
\end{pmatrix}
\end{aligned}
\end{equation}
Based on the redefined rate of the incremental deformation gradient $\Delta \dot{\tenF}^i$ (\ref{eq:rate_Berg}), the incremental deformation gradient (\ref{eq:defo_grad}), the velocity divergence (\ref{eq:div_v}) and the gradient of the pressure (\ref{eq:grad_p}) are computed. The loss function reduces to the contributions from the IRK integration of the velocity (\ref{eq:loss}), the incompressibility constraint (\ref{eq:loss_div}) and the pressure Dirichlet boundary conditions (\ref{eq:loss_boun})$_2$:
\begin{equation}
\begin{aligned}
	\mathcal{L} = SSE_{\vecv} + SSE_{\operatorname{div}\tenv} + SSE_{\bar{p}}
\end{aligned}
\end{equation}
The static pressure field obtained with the boundary projection is shown in figure \ref{fig:BC_training}. {Using the developed approach, the exact pressure is accurately computed even on a highly irregular spatial discretization.} This is a significant advantage over traditional mesh-free methods.

\subsection{Contact with rigid walls}

The velocity projection introduced in the preceding section ensures the fulfillment of Dirichlet boundary conditions for particles that belong to the boundary. However, it does not prevent fluid particles from crossing the boundary in the case of large deformations as occurring in the dam break example of section \ref{subsec:damBreak}. An additional algorithm to treat the contact of the fluid with a rigid wall is required. We choose the penalty approach which can be regarded as an artificial volumetric spring force $\tenf_c^i$. The elongation of the spring is the signed gap $g^i$ between the prescribed boundary position $\bar{\tenx}$ and the predicted current IRK stages of the position $\tenx^i$. It is penalized by the parameter $\varepsilon_c$:
\begin{equation}\label{eq:penalty_force}
\begin{aligned}
	\tenf_c^i = \varepsilon_c \, g^i \, a\left(g^i \right) \vecn, 
	\qquad 
	g^i = \left(\tenx^i-\bar{\vecx}\right)\cdot \vecn
\end{aligned}
\end{equation}
The vector $\tenn$ is the outer surface normal. The spring is only active, when a particle has crossed a boundary. This is modeled by an activation function $a\left(g^i\right)$ which takes values in the range from 0 (no contact) to 1 (contact):
\begin{equation}
\begin{aligned}
	a\left(g^i\right) = 0.5\left(1 + \operatorname{sign}\left(g^i\right)\right) = \begin{cases}
	0   & \text{if } g^i< 0 \\
	0.5 & \text{if } g^i = 0 \\
	1   & \text{if } g^i > 0
	\end{cases}
\end{aligned}
\end{equation}
For further information on the penalty method {and contact mechanics}, the reader is referred to standard text books, e.g. \cite{Wriggers.2008}.

\section{Numerical Results}\label{sec:results}

The performance of the Neural Particle Method is illustrated by means of numerical examples. To demonstrate its excellent conservation properties, free sloshing oscillations in a container are simulated for the inviscid case. Both small and large amplitude sloshing can be displayed. The classical dam break test case is simulated until the fluid hits the opposite wall. The simulations are in excellent agreement with experimental results reported in the literature.

\subsection{Inviscid free sloshing}\label{subsec:sloshing}

\begin{figure}
	\centering
	\includegraphics[width=0.55\textwidth]{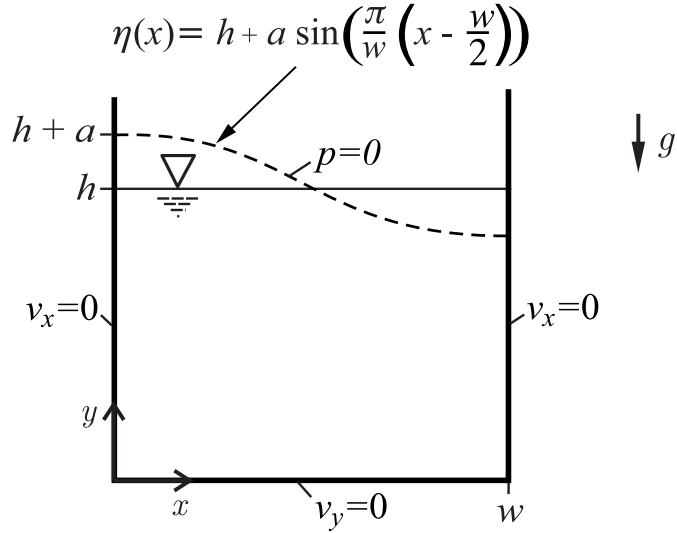}
	\caption{Setup of the sloshing test case. In the present work, the geometry is defined by $w=h=1$ m and $a=0.01$ m. The density of the liquid is $\rho=1$ kg/m$^3$ and the gravity acceleration $g=1$ m/s$^2$.}
	\label{fig:setup_sloshing}
\end{figure}

\begin{figure}
	\centering
%	\subfigure{\includegraphics[width=0.49\textwidth]{./figures/sloshing/paper/analytical_1e-1.png}}
	\subfigure{\includegraphics[width=0.9\textwidth]{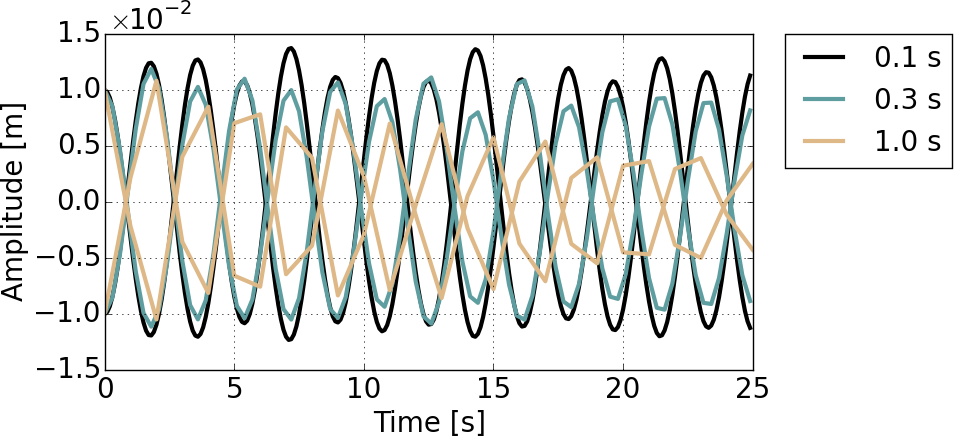}}
	
	\caption{\textbf{Constant network layout - variable time step}: Amplitude-time diagram for the inviscid sloshing test case with different time steps. The network layout used in all simulations consists of 2 input neurons, 2 hidden layers with 60 neurons each and 62 output neurons, i.e. $q=20$ IRK stages. This corresponds to the layout 1 from table \ref{tab:network_layouts}.}
	\label{fig:amplitude_inviscid}
	
	\subfigure{\includegraphics[width=0.49\textwidth]{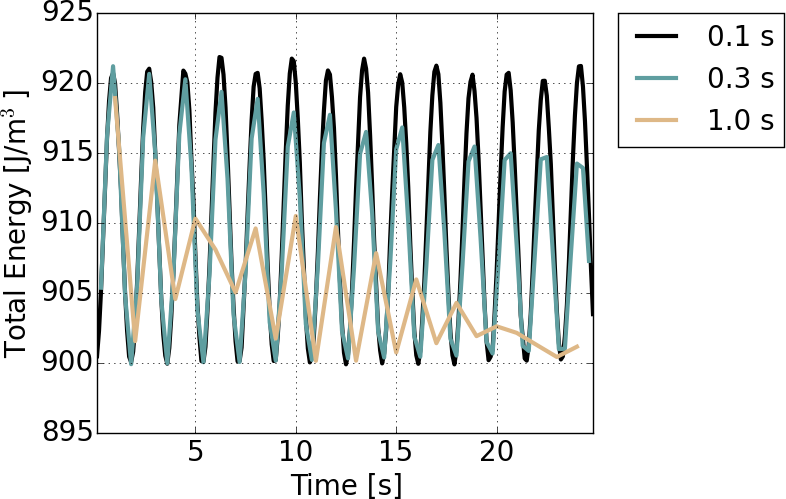}}
	\subfigure{\includegraphics[width=0.49\textwidth]{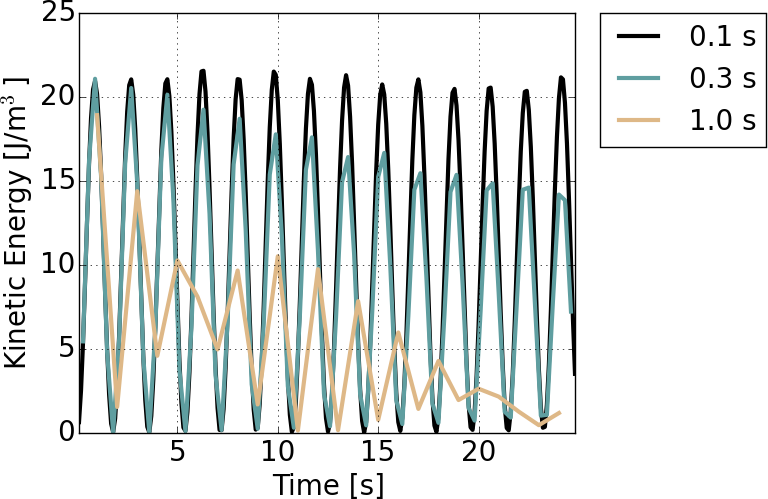}}
	
	\caption{\textbf{Constant network layout - variable time step}: Energy-time diagram (left) and kinetic energy-time diagram (right) for the inviscid sloshing test case with different time steps. The network layout is the same as in figure \ref{fig:amplitude_inviscid}.}
	\label{fig:energy_inviscid}

	\subfigure{\includegraphics[width=0.49\textwidth]{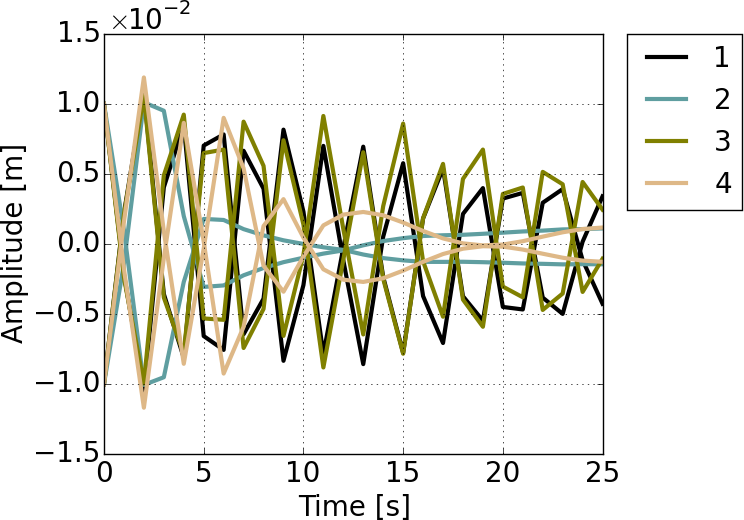}}
	\subfigure{\includegraphics[width=0.49\textwidth]{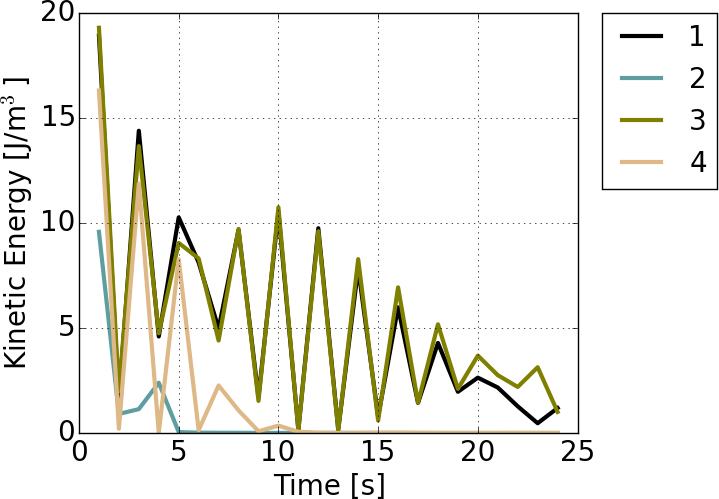}}
	
	\caption{\textbf{Constant time step - variable network layout}: Amplitude-time diagram (left) and kinetic energy-time diagram (right) for the inviscid sloshing test case with different number of IRK stages for a large time step $\Delta t=1$s. The labels correspond to different network layouts summarized in table \ref{tab:network_layouts}.}
	\label{fig:IRK_stages_inviscid}
\end{figure}

Sloshing in an open container of unit width $w$ and unit height $h$ is considered. The fluid has unit density and is subject to gravity of a unit magnitude. Initially, the surface elevation $\eta$ of the fluid follows a sine profile with amplitude $a$:
\begin{equation}\label{eq:surface_elevation}
\begin{aligned}
	\eta\left(x\right) = h - a \, \sin{\frac{\pi}{w}\left( x - \frac{w}{2}\right)}
\end{aligned}
\end{equation}
The normal velocity to the side and bottom walls is equal to zero as well as the pressure on the fluid free surface. The geometry and the boundary conditions are graphically summarized in figure \ref{fig:setup_sloshing}. This test case is especially interesting because for small sloshing amplitudes an analytical solution exists. It is widely used as benchmark test case in the literature, see e.g. \cite{Ramaswamy1990}, \cite{RadoOrtiz1998}, \cite{Braess2000} or \cite{Onate2004}. In order to evaluate the conservation properties of the method, {the total energy of the system is observed}. The {specific total energy} is comprised of a pressure, a kinetic and a potential contribution: 
\begin{equation}\label{eq:energy}
\begin{aligned}	
E^{tot} = p + \onetwo \rho \, \vecv^2 +  \rho \, g \, h
\end{aligned}
\end{equation}
In order to study the conservation properties of NPM, small amplitudes $a = 0.01$ m are considered first. Unless otherwise stated, the spatial discretization consists of 900 nearly equi-spaced points where the height coordinate of all inner and surface particles is shifted according to the surface elevation (\ref{eq:surface_elevation}). \newline

\textbf{Constant network layout - variable time step}\newline

To investigate the convergence of the method with respect to the temporal discretization, the network layout is kept constant. It consists of 2 input neurons, 2 hidden layers with 60 neurons each and 62 output neurons, i.e. $s=20$ IRK stages. The first 14 sloshing amplitudes for the inviscid case are plotted in figure \ref{fig:amplitude_inviscid} for different time step sizes $\Delta t \in \left[0.1, \, 1\right]$ s. When looking at the energy-time diagram in figure \ref{fig:energy_inviscid}, it is observed that the energy is only conserved for the small time step $\Delta t=0.1$ s. For a large time step $\Delta t = 1$ s, artificial damping becomes dominant and the amplitude decays rapidly. Whilst the variation in pressure and potential energy is below $1\%$ in all cases, the kinetic energy is significantly decaying. Next, it is examined whether the restrictions on the time step size can be relaxed if more hidden neurons are added to the network. \newline

\textbf{Constant time step - variable network layout}\newline

Different network architectures summarized in table \ref{tab:network_layouts} have been examined in combination with a large time step $\Delta t=1$ s. The computed amplitude as well as the kinetic energy of the system are plotted in figure \ref{fig:IRK_stages_inviscid}. It is found that increasing the number of IRK stages and the complexity of the neural network does not overcome the restrictions on the time step size. This result is in sharp contrast to the the findings of \cite{Raissi2019}, who presented the IRK integration method as a versatile tool enabling the use of arbitrarily large time steps.
\begin{figure}
	\centering
	
	\subfigure{\includegraphics[width=0.22\textwidth]{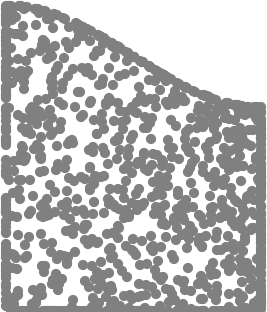}}
	\hspace{2mm}
	\subfigure{\includegraphics[width=0.22\textwidth]{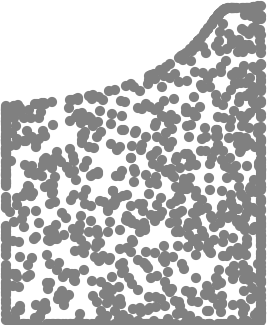}}
	\hspace{2mm}
	\subfigure{\includegraphics[width=0.22\textwidth]{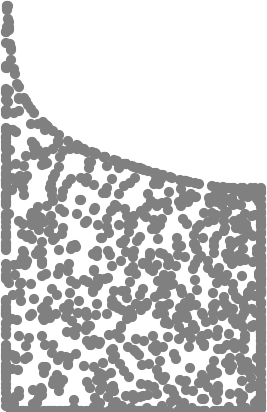}}
	\hspace{2mm}
	\subfigure{\includegraphics[width=0.22\textwidth]{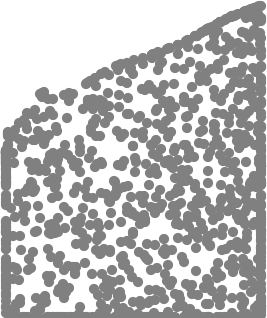}}
	\hspace{2mm}
	
	\caption{Inviscid sloshing with high amplitude $a=0.2$ m at T = 0 s, 1.9 s, 3.7 s, 5.7 s (from left to right).}
	\label{fig:sloshing_high_ampl}
\end{figure}
\begin{table}[b]
	\centering
	\begin{tabular}{lll}
		Label	& Layout& IRK stages  \\
		\hline
		1	& $\left[2, 60, 60, 62\right]$ &20  \\
		2	& $\left[2, 100, 100, 100, 100, 100, 100, 152\right]$ &50  \\
		3   & $\left[2, 200, 200, 200, 152\right]$ &50 \\
		4	& $\left[2, 300, 300, 300, 300, 300, 300, 152\right]$ & 50
	\end{tabular}
	\caption{List of network layouts used for the results presented in figure \ref{fig:IRK_stages_inviscid}.}
	\label{tab:network_layouts}
\end{table}
This limitation is likely due to the underlying assumption of the Updated Lagrangian formulation on the spatial derivatives, namely the push-forward operation using the incremental deformation gradient:
\begin{equation}
\begin{aligned}
	\frac{\partial \bullet}{\partial \vecx\na} \approx \frac{\partial \bullet}{\partial \vecx_n} \Delta \tenF^{-1},
	\qquad
	\Delta \tenF = \frac{\partial \vecx_n}{\partial \vecx\na}
\end{aligned}
\end{equation}
Note that the deformation gradient is a linear operator. It can be derived from a Taylor series expansion of the displacement, see e.g. \cite{Spencer2004}. The current position can be expressed in terms of the position at time $t_n$ and the displacement increment between both configurations as $\Delta \tenu\na = \vecu\na - \vecu_n$. In differential form this yields:
\begin{equation}
\begin{aligned}
	\mbox{d}\tenx\na &= \mbox{d}\tenx_n + \mbox{d}\Delta \tenu\na \\
	&= \left(\tenone + \frac{\partial \Delta \tenu\left(\tenx_n\right)}{\partial \tenx_n}\right) \mbox{d}\tenx_n + \frac{\partial^2 \tenu\left(\tenx_n\right)}{\partial \tenx_n^2} \, \mbox{d}\tenx_n\otimes \mbox{d}\tenx_n + ...\\
	&= \Delta \tenF \, \mbox{d}\tenx_n + \frac{\partial^2 \Delta \tenu\left(\tenx_n\right)}{\partial \tenx_n^2} \, \mbox{d}\tenx_n \otimes \mbox{d}\tenx_n + ...
\end{aligned}
\end{equation}
While generally higher order approximations are possible, these necessitate a meaningful definition of the line increment $\mbox{d}\tenx_n$ which is beyond the scope of the present work.\newline

\textbf{Sloshing at large amplitudes}\newline

The sloshing amplitude has been increased about a factor of twenty to $a=0.2$ m. The spatial discretization consists of 100 fluid particles at each fluid boundary and 700 interior particles. All points are randomly distributed. The network layout 1 from table \ref{tab:network_layouts} and a time step $\Delta t = 0.1$ s are employed. For such a large amplitude, additional waves are overlapping causing the wave to break.  The numerical results are displayed in figure \ref{fig:sloshing_high_ampl}. They  are in good agreement with the {results} reported by \cite{Onate2004} obtained with the Particle Finite Element Method.

\begin{figure}[t]
	\centering
	\subfigure{\includegraphics[width=0.4\textwidth]{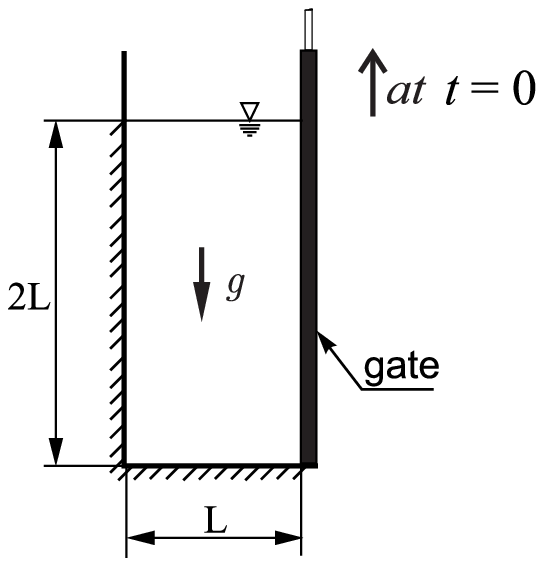}}
	
	\caption{Dam break simulation setup adapted from \cite{Ramaswamy1987}.}
	\label{fig:damBreak_setup}
\end{figure}

\begin{figure}
	\centering
	
	\subfigure{\includegraphics[width=0.7\textwidth]{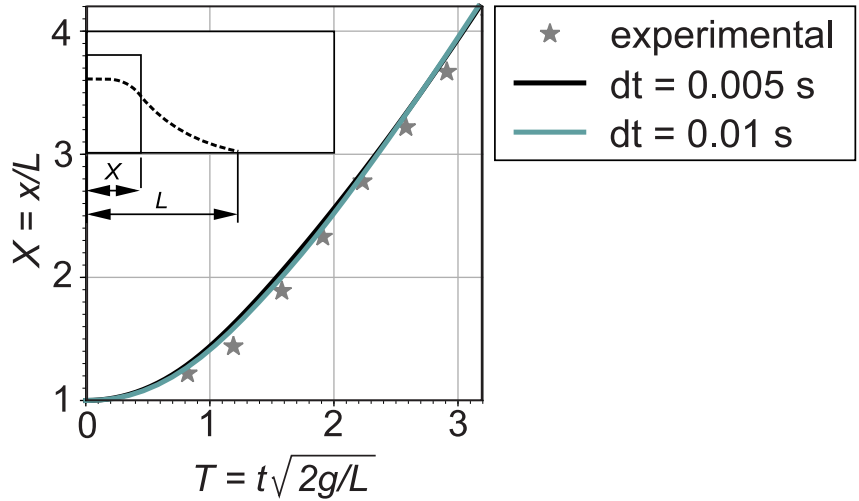}}
	
	\caption{Comparison of measured and simulated evolution of the water front tip. The experimental data is taken from \cite{Martin1952}. The training set consists of 20 particles per unit $L$. No significant difference in results could be observed to the same simulations with 25 particles per unit $L$.}
	\label{fig:experimental}
\end{figure}

\begin{figure}
	\centering
	
	\subfigure{\includegraphics[width=0.45\textwidth]{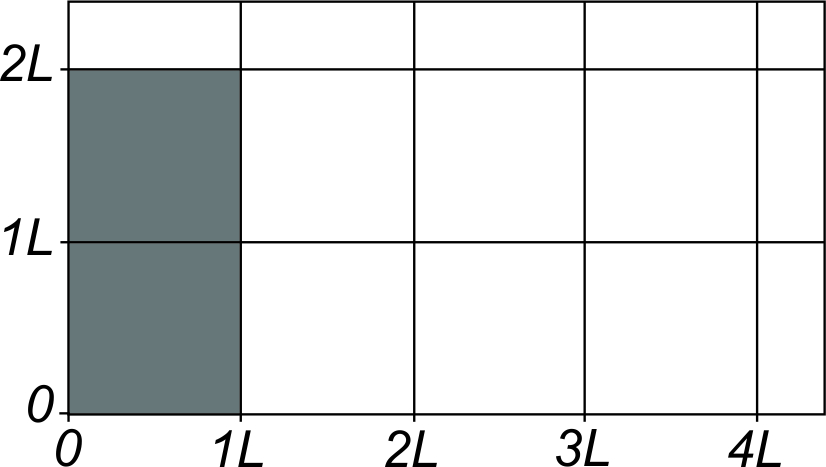}}
	\hspace{7mm}
	\subfigure{\includegraphics[width=0.45\textwidth]{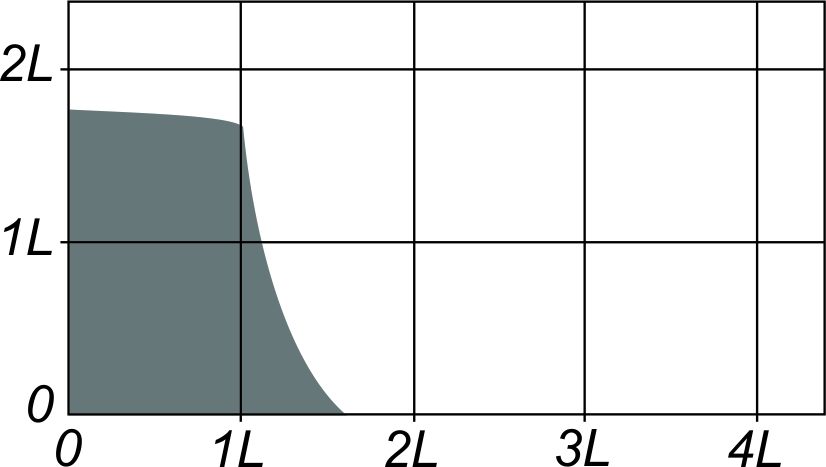}}
	
	\vspace{7mm}
	
	\subfigure{\includegraphics[width=0.45\textwidth]{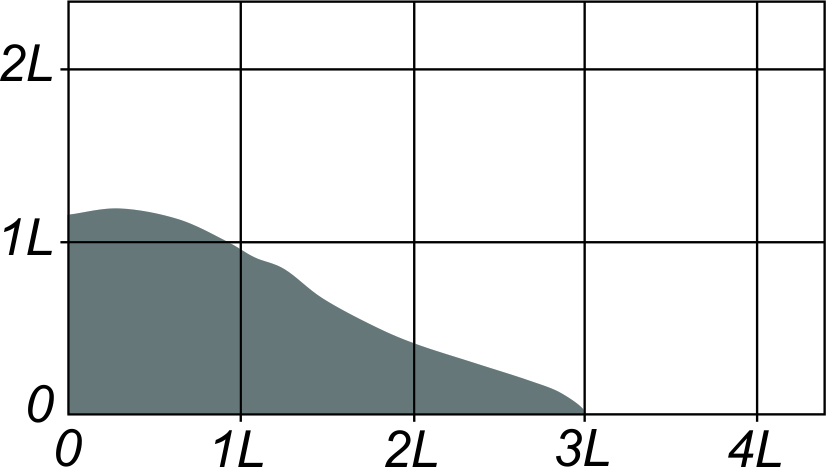}}
	\hspace{7mm}
	\subfigure{\includegraphics[width=0.45\textwidth]{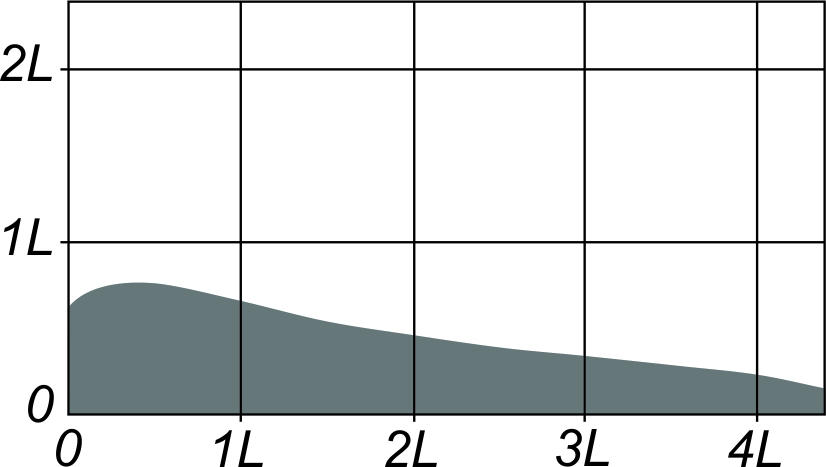}}
	
	\caption{Simulation results at T = 0 s, 0.1 s, 0.2 s, 0.3 s (left to right). The results agree well with those reported by \cite{Koshizuka1996}. The time step is $\Delta t = 0.01$ and the discretization consists of 25 particles per unit $L$. Only a subset of the displayed points was used for training.}
	\label{fig:damBreak_shape}
\end{figure}

\begin{table}[t]
	\centering
	\begin{tabular}{lll}
		Parameter	& Symbol & Value \\
		\hline
		Gravity acceleration & $g$ 				& 9.8 m/s$^2$  \\
		Length				 & $L$ 				& 0.146 m  \\
		Penalty parameter 	 & $\varepsilon_c$  & $10^7$ \\
	\end{tabular}
	\caption{Parameter used in the dam break problem.}
	\label{tab:network_layouts}
\end{table}

\subsection{Dam break}\label{subsec:damBreak}

The dam break example is a classical validation test case for Lagrangian fluid simulations, see e.g. \cite{Hirt.1974}, \cite{Ramaswamy1987} or \cite{Koshizuka1996}. Tabulated experimental data has been reported by \cite{Martin1952}. The geometry is sketched in figure \ref{fig:damBreak_setup}. Note that in this paper only the outflow of the water column has been simulated. In order to simulate the fluids interaction with an opposite wall, an additional surface identification algorithm would be required to ensure a correct imposition of the zero pressure boundary condition. This may be realized for example with the alpha-shape technique which is based upon a Delaunay triangulation of the domain, see \cite{Edelsbrunner1994}. However, surface identification is beyond the scope of the present work and only the outflow of the fluid is considered.

Beside the pressure boundary, special attention must be paid at the velocity boundaries to ensure a proper imposition of the slip condition. In SPH, dummy particles with zero velocity are used to enforce the boundary constraints. With a neural network as global ansatz, this boundary treatment corresponds to a stick condition and is not feasible. But also the velocity projection (\ref{eq:velo_proj}) would result in zero velocity at the container edges and result in a stick condition. To allow the inviscid fluid to slip along the container walls, the velocity projection is relaxed. If the distance of a fluid particle to an edge falls below a threshold, the particle is shifted along the wall in direction of the fluid motion. It then becomes part of the boundary which is perpendicular to the original one. Note that a linear velocity projection (\ref{eq:velo_proj})$_2$ is used in both spatial directions, where the width and the height of the domain are computed prior to each time step.

Figure \ref{fig:experimental} demonstrates an excellent agreement of the simulated fluid front tip evolution with the experimental data reported by \cite{Martin1952}. Network layout 1 from table \ref{tab:network_layouts} has been employed in all simulations.  The results have been obtained on a equispaced discretization of 20 particles per unit $L$. A refinement to 25 particles per unit $L$ did not affect the results. The overall shape of the collapsed water column in the time interval $T=\left[0, 0.3\right]$ s plotted in figure \ref{fig:damBreak_shape} also agrees well with the simulation results of \cite{Koshizuka1996}. Only difference is the reduced wettability of the fluid near the vertical wall. This may be due to the imposition of the zero pressure boundary condition on the free surface, which is not updated in the absence of a surface identification algorithm as mentioned above.

\section{Conclusion}\label{sec:conclusion}

The Neural Particle Method (NPM) has been proposed as a versatile simulation tool for incompressible fluid flow involving free surfaces. A feed forward neural network is chosen as spatial ansatz function for the velocity and the pressure \cite[]{Lagaris1998}. Boundary conditions are exactly fulfilled due to {the implementation of a boundary projection method} introduced by \cite{Berg2018}. For the temporal integration, high order Implicit Runge Kutta (IRK) methods are {applied} \cite[]{Raissi2019}. The inviscid, incompressible Euler equations were formulated in an Updated Lagrangian manner. NPM {computes} the pressure that {accurately} fulfills the incompressibility constraint while topological restrictions on the discretization are not required. This superior behavior is in sharp contrast to state of the art numerical methods which easily fail when it comes to disordered particle configurations. The conservation properties of the method were demonstrated in a sloshing test case. Additionally it was shown {that} the method performs well even for large sloshing amplitudes.

Although high order IRK methods were used, the admissible time step size {for the problems at hand is limited}. In the Updated Lagrangian framework, spatial derivatives are pushed into the current configuration with an incremental deformation gradient. However, the latter is only a linear operator which restricts the time step size. In order to exploit the full potential of the IRK integration, future work should focus on a relaxation of the time step constraint originating from the linear deformation map.

\section{Data availability}

The code for all numerical examples can be downloaded from https://gitlab.com/henningwessels/npm.

\section{Acknowledgment}

The first author wants to thank Jan Niklas Fuhg for the fruitful discussions.

\appendix

\bibliographystyle{plainnat}
\bibliography{npm}

%End of document
\end{document}